\def\LaTeX{%
  \let\Begin\begin
  \let\End\end
  \let\salta\relax
  \let\finqui\relax
  \let\futuro\relax}

\def\UK{\def\our{our}\let\sz s}
\def\USA{\def\our{or}\let\sz z}



\LaTeX

\USA


\salta

\documentclass[12pt]{article}  
\setlength{\textheight}{23.5cm}
\setlength{\textwidth}{16cm}
\setlength{\oddsidemargin}{2mm}
\setlength{\evensidemargin}{2mm}
\setlength{\topmargin}{-8mm}
\parskip2mm


\usepackage{color}
\usepackage{amsmath}
\usepackage{amsthm}
\usepackage{amssymb}
\usepackage[mathcal]{euscript}





\bibliographystyle{plain}


%

\finqui

\def\Beq{\Begin{equation}}
\def\Eeq{\End{equation}}
\def\Bsist{\Begin{eqnarray}}
\def\Esist{\End{eqnarray}}

\def\Bthm{\Begin{theorem}}
\def\Ethm{\End{theorem}}
\def\Blem{\Begin{lemma}}
\def\Elem{\End{lemma}}

\def\Brem{\Begin{remark}\rm}
\def\Erem{\End{remark}}

\def\Bdim{\Begin{proof}}
\def\Edim{\End{proof}}
\let\non\nonumber




\def\step #1 \par{\medskip\noindent{\bf #1.}\quad}


\def\aand{\quad\hbox{and}\quad}
\def\where{\quad\hbox{where}\quad}

\def\lhs{left-hand side}
\def\rhs{right-hand side}
\def\sfw{straightforward}



\def\multibold #1{\def\arg{#1}%
  \ifx\arg\pto \let\next\relax
  \else
  \def\next{\expandafter
    \def\csname #1#1#1\endcsname{{\bf #1}}%
    \multibold}%
  \fi \next}

\def\pto{.}

\def\multical #1{\def\arg{#1}%
  \ifx\arg\pto \let\next\relax
  \else
  \def\next{\expandafter
    \def\csname #1#1\endcsname{{\cal #1}}%
    \multical}%
  \fi \next}

\def\multizero #1{\def\arg{#1}%
  \ifx\arg\pto \let\next\relax
  \else
  \def\next{\expandafter
    \def\csname #1z\endcsname{#1_0}%
    \multizero}%
  \fi \next}

\def\multistar #1{\def\arg{#1}%
  \ifx\arg\pto \let\next\relax
  \else
  \def\next{\expandafter
    \def\csname #1star\endcsname{#1^*}%
    \multistar}%
  \fi \next}

\let\hat\widehat
\let\tilde\widetilde

\def\multihat #1{\def\arg{#1}%
  \ifx\arg\pto \let\next\relax
  \else
  \def\next{\expandafter
    \def\csname hat#1\endcsname{\hat #1}%
    \multihat}%
  \fi \next}


\def\multimathop #1 {\def\arg{#1}%
  \ifx\arg\pto \let\next\relax
  \else
  \def\next{\expandafter
    \def\csname #1\endcsname{\mathop{\rm #1}\nolimits}%
    \multimathop}%
  \fi \next}

\multibold
qwertyuiopasdfghjklzxcvbnmQWERTYUIOPASDFGHJKLZXCVBNM.

\multical
QWERTYUIOPASDFGHJKLZXCVBNM.

\multizero
qwertyuiopadfghjklzxcvbnmQWERTYUIOPASDFGHJKLZXCVBNM.  

\multistar
qwertyuiopasdfghjklzxcvbnmQWERTYUIOPASDFGHJKLZXCVBNM.

\multihat
qwertyuiopasdfghjklzxcvbnmQWERTYUIOPASDFGHJKLZXCVBNM.

\multimathop
dist div dom meas sign supp .


\def\accorpa #1#2{\eqref{#1}--\eqref{#2}}
\def\Accorpa #1#2 #3 {\gdef #1{\eqref{#2}--\eqref{#3}}%
  \wlog{}\wlog{\string #1 -> #2 - #3}\wlog{}}


\def\graffe #1{\mathopen\{#1\mathclose\}}

\def\<#1>{\mathopen\langle #1\mathclose\rangle}
\def\norma #1{\mathopen \| #1\mathclose \|}
\def\modulo #1{\left| #1\right|}

\def\ioT {\int_0^T}

\def\dt{\partial_t}
\def\dn{\partial_n}

\def\checkmmode #1{\relax\ifmmode\hbox{#1}\else{#1}\fi}
\def\aeO{\checkmmode{a.e.\ in~$\Omega$}}
\def\aeQ{\checkmmode{a.e.\ in~$Q$}}

\def\aaO{\checkmmode{for a.a.~$x\in\Omega$}}
\def\aaQ{\checkmmode{for a.a.~$(x,t)\in Q$}}


\def\erre{{\mathbb{R}}}




\def\genspazio #1#2#3#4#5{#1^{#2}(#5,#4;#3)}
\def\spazio #1#2#3{\genspazio {#1}{#2}{#3}T0}

\def\spazioi #1#2#3{\genspazio {#1}{#2}{#3}\infty 0}
\def\L {\spazio L}
\def\H {\spazio H}
\def\W {\spazio W}

\def\LL {\spazioi L}

\def\C #1#2{C^{#1}([0,T];#2)}


\def\Lx #1{L^{#1}(\Omega)}
\def\Hx #1{H^{#1}(\Omega)}
\def\Wx #1{W^{#1}(\Omega)}
\def\Luno{\Lx 1}
\def\Ldue{\Lx 2}
\def\Linfty{\Lx\infty}
\def\Huno{\Hx 1}


\def\LinftyQ{L^\infty(Q)}
\def\LdueQ{L^2(Q)}
\def\LunoQ{L^1(Q)}


\let\theta\vartheta
\let\eps\varepsilon
\let\phi\varphi

\let\TeXchi\chi                         
\newbox\chibox
\setbox0 \hbox{\mathsurround0pt $\TeXchi$}
\setbox\chibox \hbox{\raise\dp0 \box 0 }
\def\chi{\copy\chibox}


\let\u\rho
\def\uz{\u_0}

\def\xiz{\xi_0}
\def\ximax{{\xi^*}}

\let\emme\sigma

\def\epsz{\eps_0}
\def\deltaz{\delta_0}
\def\Lz{L_0}
\def\etaz{\eta_0}

\def\rhoz{\rho_0}
\def\rhomin{{\rho_*}}
\def\rhomax{{\rho^*}}
\def\sqrho{\sqrt\rho}
\def\sqrhoxi{\sqrt{\rho\,\xi}}
\def\sqrhonxin{\sqrt{\rhon\,\xin}}
\def\rhon{\rho_n}

\def\xiz{\xi_0}
\def\ximax{{\xi^*}}

\def\xin{\xi_n}

\def\thetaz{\theta_0}
\def\thetam{\underline s}
\def\thetaM{\bar s}
\def\thetamin{\theta_{\!*}}
\def\thetamax{\theta^*}
\def\thetan{\theta_{\!n}}
\def\htheta{\bar\theta}
\def\hthetan{\bar\theta_{\!n}}

\def\source{\bar\sigma}

\def\fu{f_1}
\def\fd{f_2}

\def\lavoro{\cite{CGPS}}


\def\cv{c_v}
\def\cs{c_*}
\DeclareMathAlphabet{\mathbf}{OT1}{cmr}{bx}{it}

\newcommand{\0}{\mathbf 0}
\newcommand{\hb}{\mathbf{h}}
 \font\mba=cmmib10 scaled
\magstephalf

\def\csib{\hbox{\mba {\char 24}}}
\newcommand{\pier}{}

\topmargin=0cm
\oddsidemargin=0cm
\textwidth=16cm
\textheight=23cm
\topmargin=-1cm 
\Begin{document}

\begin{center}
{\Large\bf A temperature-dependent phase segregation problem\\[0.25cm] of the Allen-Cahn type}\\[1cm] 
{\large\bf Pierluigi Colli$^{(1)}$}\\
{\normalsize e-mail: {\tt pierluigi.colli@unipv.it}}\\[.4cm]
{\large\bf Gianni Gilardi$^{(1)}$}\\
{\normalsize e-mail: {\tt gianni.gilardi@unipv.it}}\\[.4cm]
{\large\bf Paolo Podio-Guidugli$^{(2)}$}\\
{\normalsize e-mail: {\tt ppg@uniroma2.it}}\\[.4cm]
{\large\bf J\"urgen Sprekels$^{(3)}$}\\
{\normalsize e-mail: {\tt sprekels@wias-berlin.de}}\\[1cm]
$^{(1)}$
{\small Dipartimento di Matematica ``F. Casorati'', Universit\`a di Pavia}\\
{\small via Ferrata 1, 27100 Pavia, Italy}\\[.2cm]
$^{(2)}$
{\small Dipartimento di Ingegneria Civile, Universit\`a di Roma ``Tor Vergata''}\\
{\small via del Politecnico 1, 00133 Roma, Italy}\\[.2cm]
$^{(3)}$
{\small Weierstra\ss-Institut f\"ur Angewandte Analysis und Stochastik}\\
{\small Mohrenstra\ss e\ 39, 10117 Berlin, Germany}\\[.8cm]
{\large\it Dedicated to Professor Nobuyuki Kenmochi}\\[.2cm]
{\large\it on the occasion of his 65th birthday}
\end{center}

\bigskip
\begin{abstract}
In this paper we prove a local-in-time existence theorem for an initial-boundary 
value problem related to a model of temperature-dependent phase segregation 
that generalizes the standard Allen-Cahn's model. The problem is ruled 
by a system of two differential equations, one partial the other ordinary,
interpreted as balances, respectively, of microforces and of microenergy,
complemented by a transcendental condition on the three
unknowns, that are: the order parameter entering the standard A-C equation, the
chemical potential, and the absolute temperature. The results obtained 
in our recent paper \lavoro\  dealing with the isothermal case serve as a 
starting point for our existence proof, which relies on a
fixed-point argument involving the Tychonoff-Schauder theorem.
\end{abstract}





\section{Introduction}
\label{Intro}
\setcounter{equation}{0}
 This paper is a sequel and a generalization of our article \lavoro, where 
we studied a nonlinear evolution system of the Allen-Cahn (A-C) type,  
intended to provide a mathematical description of the phenomenology of 
phase segregation by atom rearrangement in crystalline materials, in the 
absence of diffusion. Our present generalization consists in taking thermal 
effects into account. To help understanding the underlying physics, 
temperature-independent mathematical models of A-C type and their proposed 
generalization are briefly discussed in the next section, where we recapitulate 
Gurtin's derivation of the standard A-C equation as well as the derivation of 
the nonstandard A-C system analyzed in \lavoro, and where we sketch the main 
traits of our present temperature-dependent  model. In Section 3, we formulate 
carefully the corresponding mathematical problem and we state a local-in-time 
existence result, that we prove in our last section. 

\section{Phase segregation models of A-C type}
\subsection{The temperature-independent model of \lavoro}
In \lavoro, we consider a nonlinear evolution system consisting of
the partial differential equation: 
\begin{equation}  \label{Iprima}
\kappa \, \dt\u - \Delta\u + f'(\u) = \mu 
\end{equation}
and of the ordinary differential equation:  
\begin{equation} \label{secondanew}
 \dt(-\mu^2\rho) = \mu\left(\kappa \, (\dt\u)^2+\bar\emme \right) ,
\end{equation}
complemented with the homogeneous Neumann boundary condition:
\Beq   \label{Ibcic1}
  \dn \u = 0 \quad \hbox{on the body's boundary} 
  \Eeq
(here $\dn$ denotes the outward normal derivative) and with the initial conditions: 
\Beq \label{Ibcic2}
   \u|_{t=0} = \uz\quad\hbox{bounded away from} \;0 \,,\quad\mu|_{t=0} = \mu_0\geq 0\,.
\Eeq

The parabolic PDE (\ref{Iprima}) and the first-order-in-time ODE (\ref{secondanew}) are interpreted as balances, respectively, of microforces and of microenergy. They are to be solved for the order-parameter field $\rho= \rho(x,t)\in [0,1]$, interpreted as the scaled volumetric density 
of one of two coexisting phases, and for the chemical potential field $\mu$. 
Moreover,  $\kappa>0$~is a mobility coefficient, $f$ is 
a double-well potential confined in~$(0,1)$ and singular at endpoints, 
and $\bar\emme = \bar\emme (x,t)$ denotes a given source term. 
The microentropy field $\eta = -\mu^2 \rho$ cannot exceed the level 
$0$ from below, so that the corresponding prescribed initial field 
\Beq\label{eta0}
\eta|_{t=0} =\eta_0\, ,
\Eeq
with $\eta_0=-\mu_0^2\rho_0$, is nonpositive-valued.  
Taking $\mu\equiv 0$ in \eqref{Iprima} yields the standard 
Allen-Cahn equation,  which is intended to describe evolutionary processes in a two-phase 
material body, \emph{phase segregation} included. 

\subsection{Gurtin's derivation of the A-C equation}
The derivation of the 
A-C equation proposed by Gurtin~\cite{gurtin} (see also \cite{fremond} and  
\cite{miranv} for similar derivations and for discussions of related models) is based on a \emph{balance 
of contact and distance microforces}:
\begin{equation}\label{balance}
\div\csib+\pi+\gamma=0
\end{equation}
and on a dissipation inequality restricting the free-energy
growth:
\begin{equation}\label{dissipation}
\partial_t\psi\leq w,\quad w:=-\pi\,\partial_t\rho+\csib\cdot\nabla(\partial_t\rho),
\end{equation}
where the distance microforce is split in an internal part $\pi$ and an external part 
$\gamma$, the contact microforce is specified by the microscopic stress vector $\csib$, and $w$ is the 
(distance plus contact) internal microworking.  
Requesting the Coleman-Noll (C-N) compatibility \cite{CN} of  the constitutive choices: 
\begin{equation}\label{constitutive}
\pi=\hat\pi(\rho,\nabla\rho,\partial_t\rho), \quad\csib=\hat\csib(\rho,\nabla\rho,\partial_t\rho), \quad \textrm{and}\quad 
\psi=\hat\psi(\rho,\nabla\rho)=f(\rho)+\frac{1}{2}|\nabla\rho|^2\,,
\end{equation}
%
%
with the dissipation inequality \eqref{dissipation} yields:
 \[
\hat\pi(\rho,\nabla\rho,\partial_t\rho)=-f'(\rho)-\hat\kappa(\rho,\nabla\rho,\partial_t\rho)\partial_t\rho\quad\textrm{and}\quad \hat\csib(\rho,\nabla\rho,\partial_t\rho)=\nabla\rho.
\]
Under the further assumptions that  $\hat\kappa(\rho,\nabla\rho,\partial_t\rho)=\kappa$, a positive constant, and that $\gamma\equiv 0$, the microforce balance (\ref{balance}) takes the form of the standard Allen-Cahn equation
\begin{equation}  \label{ac}
\kappa \, \dt\u - \Delta\u + f'(\u) = 0.
\end{equation}

\subsection{The derivation of the A-C system of \lavoro}\label{AC}
In \lavoro, on adopting an approach put forward by one of us in~\cite{Podio}, 
we deal with a modified version of Gurtin's derivation, in which the dissipation
inequality \eqref{dissipation} is dropped and the microforce balance 
\eqref{balance} is coupled both with the \emph{microenergy balance}
\begin{equation}\label{energy}
\partial_t\varepsilon=e+w,\quad e:=-\div{\bar \hb}+{\bar \sigma},
\end{equation}
and the \emph{microentropy imbalance}
\begin{equation}\label{entropy}
\partial_t\eta\geq -\div\hb+\sigma,\quad \hb:=\mu{\bar \hb},\quad \sigma:=
\mu\,{\bar \sigma}
\end{equation}
(here $\bar\sigma$ is the external source of energy per unit volume). With a view toward modeling phase-segregation, we postulate that the microentropy inflow $(\hb,\sigma)$ be proportional to the microenergy inflow $({\bar \hb},{\bar\sigma})$ 
through the chemical potential $\mu$, a positive field; consistently, 
we define the free energy to be:
\begin{equation}\label{freeenergy}
\psi:=\varepsilon-\mu^{-1}\eta,
\end{equation}
with the chemical potential playing the same role as the coldness $\vartheta^{-1}$ 
in the deduction of the heat equation. Combination of \accorpa{energy}{freeenergy} 
leads to the inequality:
\begin{equation}\label{reduced}
\partial_t\psi\leq -\eta\,\partial_t(\mu^{-1})+\mu^{-1}{\bar \hb}\cdot\nabla\mu-\pi\,\partial_t\rho+\csib\cdot\nabla(\partial_t\rho),
\end{equation}
which replaces  \eqref{dissipation} as a restriction on constitutive choices. We assume that, in addition to the independent variables $\rho,\, \nabla\rho$ and $\partial_t \u$, the constitutive mappings delivering $\pi,\xi,\eta$, 
and ${\bar \hb}$ depend also on $\mu$; moreover, we choose:
\begin{equation}\label{constitutives}
\psi=\hat\psi(\rho,\nabla\rho,\mu)=-\mu\,\rho+f(\rho)+\frac{1}{2}|\nabla\rho|^2.
\end{equation}
To impose the C-N compatibility of these assumptions with \eqref{reduced} makes sense, because we have at our disposal two independent controls $\gamma$ and $\bar\sigma$ to guarantee the free linear continuation in time of any given process $t\mapsto (\rho,\mu)(t)$ at any fixed space point. We find:
\Bsist\label{CNc}
\hat\pi(\rho,\nabla\rho,\partial_t\rho,\mu)\displaystyle{=\mu-f'(\rho)-\hat\kappa(\rho,\nabla\rho,\partial_t\rho)\partial_t\rho,\quad \hat\csib(\rho,\nabla\rho,\partial_t\rho,\mu)=\nabla\rho,}&& \non\\[2.mm]
\phantom{\hat\pi(\rho,\nabla\rho,\partial_t\rho,\mu)}
\hat\eta(\rho,\nabla\rho,\partial_t\rho,\mu)\displaystyle{=-\mu^2\rho,\quad \hat{\bar \hb}(\rho,\nabla\rho,\partial_t\rho,\mu)\equiv \0;}\label{cn}
&&\qquad
\Esist
the last of these findings implies that the microenergy balance -- in general, a PDE -- becomes an ODE, a crucial mathematical simplification that we exploit in the following just as we did in \lavoro. Finally, under the additional assumptions that the 
mobility is a positive constant and the external distance 
microforce is null, the microforce balance \eqref{balance} 
and the energy balance \eqref{energy} become, 
respectively, \eqref{Iprima}~and~\eqref{secondanew}.

\subsection{Accounting for thermal effects}
As is well-known (see e.g. \cite{Podio}), the classic heat equation can be arrived at by coupling the energy balance
\begin{equation}\label{eb}
\partial_t\varepsilon=-\div \bar\hb
\end{equation}
and the entropy imbalance
\begin{equation}\label{eimb}
\partial_t\eta\geq-\div\hb,\quad \hb=\vartheta^{-1}\bar\hb,
\end{equation}
with the following constitutive prescriptions:
\begin{equation}\label{ce}
\psi=\varepsilon-\vartheta\eta,\quad \psi=\bar\psi(\vartheta)=-c_v\vartheta(\ln\vartheta -1),
\end{equation}
with the absolute temperature field $\vartheta$ positive-valued and the specific heat $c_v$ a positive number. To account for thermal effects on the phenomenology of phase segregation by atomic rearrangement, we compare the formats (\ref{eb})--(\ref{ce}) and (\ref{energy})--(\ref{constitutives}). A way to match them, in the light of the relationships of temperature and chemical potential to entropy provided by statistical mechanics, is  to assume that: (i) the microenergy balance keeps the form (\ref{energy}); (ii) the energy/entropy fluxes and the free energy have the mutually consistent forms
\begin{equation}\label{key}
\hb=(\vartheta^{-1}\mu)\bar\hb,\quad\psi=\varepsilon-(\vartheta\mu^{-1})\eta,\quad \psi=\widehat\psi(\rho,\nabla\rho,\mu,\vartheta).
\end{equation}
The second assumption is the main element of novelty of this note. With that measure, the dissipation inequality that replaces for (\ref{reduced}) is:
\begin{equation}\label{red}
\partial_t\psi\leq -\eta\,\partial_t(\vartheta\mu^{-1})+(\vartheta\mu^{-1}){\bar \hb}\cdot\nabla
(\theta^{-1}\mu) -\pi\,\partial_t\rho+\csib\cdot\nabla(\partial_t\rho),
\end{equation}
where, in addition to $(\ref{key})_3$ for the free energy, the
distance force, the microscopic stress, the entropy, and the microenergy
influx, are assumed to depend on the list of variables $\Lambda = \{ \rho,\nabla\rho, \mu, \theta;\partial_t \u, \nabla\theta \}$.\footnote{Note that the last two variables in the list give way to incorporate in the model the dissipation mechanisms relative to, respectively, atom-rearrangement without diffusion and heat conduction. Indeed, it is clear that (\ref{key}) covers both special cases when either
temperature or chemical potential is a space-time constant.} 

When C-N compatibility of the present constitutive prescriptions with (\ref{red}) is inspected, a delicate modeling issue emerges:  this time, we cannot count on as many controls as needed to guarantee the free local continuation in time of any given process $t\mapsto (\rho,\mu,\vartheta)(t)$ at any fixed space point. We do get the counterparts of the first, second, and fourth of (\ref{CNc}), namely,
\begin{equation}\label{rip}
\hat\pi  = - \partial_\rho \psi -\kappa \partial_t\rho,
\quad
\hat\csib =\partial_{\nabla\u} \psi = \nabla\rho,
\quad
 \hat{\bar \hb} \equiv \0;
\end{equation}
and we are left with the residual inequality:
\begin{equation}\label{redin}
(\partial_\mu\psi - \theta \mu^{-2}\eta) \partial_t \mu + 
(\partial_\theta \psi + \mu^{-1}\eta) \partial_t \theta \leq 0 .
\end{equation}
Now, if it were possible to choose both $\partial_t \mu$ and $\partial_t \vartheta$ arbitrarily, then (\ref{redin}) would yield the double equality: 
\begin{equation}\label{etaa}
\eta = \theta^{-1} \mu^2 \partial_\mu \psi = - \mu \, \partial_\theta \psi. 
\end{equation} 
This observation motivates our decision to complement the microforce and energy balances with another field equation, namely, the \emph{thermodynamic consistency condition}:
\begin{equation}\label{consi}
\mu \,\partial_\mu \psi +  \theta \, \partial_\theta \psi =0. 
\end{equation}
With this, (\ref{redin}) can be written as
\begin{equation}\label{disred}
\vartheta(\partial_\mu\psi - \theta \mu^{-2}\eta) \partial_t (\vartheta^{-1}\mu) \leq 0,
\end{equation}
and (\ref{etaa}) follows, provided the time rate of $(\vartheta^{-1}\mu)$ can be chosen arbitrarily.

Next, we specify the free energy density $(\ref{key})_3$ as follows~(cf.~\eqref{constitutives}):
\begin{equation}
\label{constitutives2}
\psi=\hat\psi(\rho,\nabla\rho,\mu, \theta)=-\mu\,\rho+ \varphi(\rho,\theta)+\frac{1}{2}|\nabla\rho|^2,
\end{equation}
with
\begin{equation}\label{fi}
 \varphi (\u, \theta ) = f(\u)  - c_v \theta( \ln \theta - 1)-  c_0\,\rho( \theta   -   \theta_c) ,\quad c_0>0,
\end{equation}
where the double-well potential and the purely caloric free energy are supplemented by a coupling term that is  effective when and where the temperature differs from the characteristic temperature  $\theta_c $. With this final constitutive choice, the consistency condition \eqref{consi} reduces to:
\begin{equation}\label{third}
\mu \rho +  c_0 \theta \u  + c_v \theta \ln \theta =0 \, ; 
\end{equation}
moreover, the balance of microforces \eqref{balance} becomes:
\begin{equation}  \label{primanew}
\kappa \, \dt\u - \Delta\u + f' (\u) - c_0 \theta = \mu
\end{equation}
where, with slight abuse of notation, we have written $f^\prime(\u)$ for $(f^\prime(\u)+c_0\theta_c)$. As to the microenergy balance (\ref{energy}), we find:
\begin{equation}\label{feb}
\partial_t(-\theta^{-1}\mu^2\rho)=\theta^{-1}\mu(\bar\sigma+\kappa(\partial_t\rho)^2),
\end{equation}
an equation to be compared with \eqref{secondanew}. Our current mathematical model regards processes of phase segregation by atomic re-arrangement in the presence of thermal effects as solutions of the system of equations (\ref{third}), (\ref{primanew}) and (\ref{feb}).

\section{Mathematical formulation and results}
\label{Results}
\setcounter{equation}{0}
The A-C system we derived is more difficult to deal with than its temperature-independent version we tackled successfully in \lavoro, let alone the standard A-C equation (\ref{ac}). In fact, in addition to a PDE and an ODE as  in \lavoro, we now have to take care also of the {\pier transcendental} equation (\ref{third}).  Our strategy is to repeat, for as much as is possible, the procedure in \lavoro: accordingly, we discuss the ODE (\ref{feb}) first together with the relative initial condition, then we pass to the PDE and the {\pier transcendental} equation, together with the relative boundary and initial conditions.

\subsection{Preliminaries}
In order to carry out the first part of
our program, we adopt a change of variable to give \eqref{feb}
plus \eqref{eta0} the form of a parametric initial-value problem. We begin by introducing 
the initial value $\theta_0 $ of $\theta$ and setting:
\[
-\eta=\theta^{-1}\xi=\theta^{-1}\mu^2\rho,\quad \xi_0=-\vartheta_0\eta_0,\quad  \eta_0= - \theta_0^{-1}\mu_0^2\rho_0.
\]
We then have that 
\Beq\label{psiz}
\mu=\sqrt{\xi/\rho},
\Eeq
whence the Cauchy problem:
\Beq\label{Cp}
\vartheta\dt(\theta^{-1}\xi)  + \frac{\kappa \, (\dt\u)^2 + \bar\emme}{\sqrt{\rho}} \, \sqrt{\xi}=0,\quad (\theta^{-1}\xi)|_{t=0} =- \eta_0.
\Eeq
Next, we restrict attention to the class  of processes such that 
\[
\vartheta\dt(\theta^{-1}\xi) \simeq \dt\xi ,
\]
and replace (\ref{Cp}) by the simpler problem:
\Beq\label{Cpp}
\dt\xi  + \frac{\kappa \, (\dt\u)^2 + \bar\emme}{\sqrt{\rho}} \, \sqrt{\xi}=0,\quad \xi|_{t=0} =\xi_0,
\Eeq
parameterized on both the space variable $x$ and the 
field $\rho(x,\cdot)$. Although simpler, this Cauchy problem is by no means trivial, because it can exhibit the 
Peano phenomenon and have infinitely many solutions; just as we did in \lavoro, we pick a suitably defined  \emph{maximal solution} $\xi $ (or $\sqrt{\xi}$), 
having the important property to stay positive as long as is possible. 

\Brem It remains to be seen whether the class we restrict attention to does include interesting phase--segregation processes. At this time, we cannot do any better than planning to check if this is the case by running numerical simulations. 
\Erem

To complete our program, we note that, with (\ref{psiz}), (\ref{primanew}) and (\ref{third}) become, respectively,
\Beq\label{acmod}
 \kappa \, \dt\u - \Delta\u + f' (\u) - c_0 \theta =   \sqrt{{\xi}/{\rho}}\,,
\Eeq
and 
\begin{equation}\label{condi}
\lambda(\rho,\theta) :=  c_0  \u  \theta + c_v \theta \ln \theta = -  \sqrt{\rho\xi}\,,
\end{equation} 
an \emph{integro--differential system} for $\rho$ and $\theta$, where $\sqrt{ \xi }$ is implicitly defined in terms of $\u$ as the maximal solution to  \eqref{Cpp}. This system is to be supplemented with the boundary condition \eqref{Ibcic1}, the initial condition for $\rho$ in  \eqref{Ibcic2}, and a compatible initial condition for~$\theta$. %
\subsection{Results}
In view of the above discussion, we look for suitably smooth triplets of time-dependent fields $(\rho,\xi,\theta)$ over a regular region $\Omega$ with boundary $\Gamma$, such that:
\Bsist
  && 0 < \rho < 1 , \quad
  \xi \geq 0 ,
  \aand
  \theta > 0;
  \label{range}
  \\
  && \dt\rho - \Delta\rho + f'(\rho) - \cz \theta
  = \sqrt{\xi/\rho}\, ,
  \quad 
  \hbox{with~ $\dn\rho=0$~ on~} \Gamma;
  \qquad
  \label{prima}
  \\
  && \dt\xi + \frac{|\dt\rho|^2 + \source}\sqrho \, \sqrt\xi = 0;
  \label{seconda}
  \\
  && \lambda(\rho,\theta) = - \sqrhoxi;
  \vrule width0pt depth10pt
  \label{terza}
  \\
  && \rho(0) = \rhoz, \quad
  \xi(0) = \xiz \,,
  \aand
  \theta(0) = \thetaz \,;
  \label{cauchy}
\Esist 
and that, moreover, 
\Beq
  \hbox{$\xi$ is maximal among the $\xi$'s
    satisfying \eqref{seconda} and the second of~\eqref{cauchy}}.
  \label{maximality}
\Eeq
\Accorpa\pbl prima maximality
The problem's structure is the same as in \lavoro, apart for the modifications due to the presence of the temperature variable $\theta$. Two items deserve a supplement of discussion. The first is that, just as in \lavoro, we assume that
\Bsist
  && 0 \leq f = \fu + \fd,
  \where
  \fu, \fd : (0,1) \to \erre
  \quad \hbox{are $C^2$-functions,}
  \label{hpfa}
  \\
  && \hbox{$\fu$ is convex}, \quad
  \hbox{$\fd'$ is bounded}, \quad
  \lim_{r\searrow 0} f'(r) = - \infty ,
  \aand
  \lim_{r\nearrow 1} f'(r) = + \infty, 
  \qquad\quad
  \label{hpfb}
\Esist
\Accorpa\structure hpfa hpfb
with the constant $\cz\theta_c$ thought of as incorporated in~$\fd'(\rho)$. 
The second item has to do with the admissible choices of initial data: not only they must agree with (\ref{range}), and hence satisfy
\Beq
  0 < \rhoz < 1 , \quad
  \xiz \geq 0 , \quad
  \thetaz > 0 , 
  \label{compdati1}
\Eeq
but also with (\ref{terza}), that is to say, they have to satisfy
\Beq
  \lambda(\rhoz,\thetaz) = - \sqrt{\rhoz\xiz} \, . 
  \label{compdati2}
\Eeq
To see what restrictions this last condition implies on the choice of $\theta_0$, it is convenient to study the function the function $\lambda_r:s\mapsto\lambda(r,s)= \cz rs + \cv s\ln s$
on $(0,+\infty)$ for a given $r\in(0,1)$.
Clearly, (i) $\lambda_r$~is strictly convex and tends to $0$ as $s$ tends to~$0$; moreover, (ii) the equation $\lambda_r(s)=0$ has in~$(0,+\infty)$  a unique solution, that we denote by $\thetaM(r)$; finally, (iii) $\lambda_r$~has a unique minimum point, denoted by $\thetam(r)$; in summary, for each fixed $r\in(0,1)$,
\Beq
  0 < \thetam(r) < \thetaM(r), \quad
  \lambda(r,\thetaM(r)) = 0 ,
  \aand
  \frac {\partial\lambda} {\partial s} (r,\thetam(r)) = 0 .
  \label{defthetamM}
\Eeq
A simple computation shows that
\Bsist
  && \thetam(r) = e^{-1-\cs r} , \quad
  \thetaM(r) = e^{-\cs r} , \aand
  \non
  \\
  && \lambda(r,\thetam(r)) = - \cv e^{-1-\cs r} ,
  \quad \hbox{where} \quad
  \cs := \cz/\cv \,.
  \label{thetamM}
\Esist
Therefore, a necessary condition for the existence of $\thetaz$
is that 
\Beq\label{necon}
\sqrt{\rhoz\xiz}\leq\cv e^{-1-\cs\rhoz}\quad \aeO,
\Eeq
 i.e., that
$\sup\zeta\leq0$, where $\zeta:=\sqrt{\rhoz\xiz}-\cv e^{-1-\cs\rhoz}$.
If such a condition is satisfied, and if we want to solve
(\aeO) the~equation $\lambda(\rhoz(x),s)=0$ for~$s$, then 
uniqueness holds if $\zeta(x)=0$,
and $\thetam(\rhoz(x))$ is the unique solution.
Otherwise, if the strict inequality holds, then there are two solutions, the one
in the interval $(0,\thetam(\rhoz(x)))$ the other in $(\thetam(\rhoz(x)),\thetaM(\rhoz(x)))$.

For  existence of a local-in-time solution~$(\rho,\xi,\theta)$, a modest reinforcement of condition (\ref{necon})
and a proper choice of  $\thetaz$ suffice, namely,
\Beq
  \sup \bigl(\sqrt{\rhoz\xiz} - \cv e^{-1-\cs\rhoz} \bigr) < 0 
  \aand
  \thetaz \geq \thetam(\rhoz)
  \quad \aeO.
  \label{hpcauchy}
\Eeq
Under these assumptions, we can state the following result.

\Bthm
\label{Existence}
Assume that \structure~and  \eqref{hpcauchy} hold.
Moreover, assume that
\Bsist
  && \source\in L^\infty(\Omega\times(0,+\infty)), \quad
  (\source)^- \in \LL\infty\Luno, 
  \aand
  \rhoz , \xiz , \thetaz \in \Linfty;
  \qquad\quad
  \label{hpdati}
  \\
  && \rhoz \in \Hx3, \quad
  \dn\rhoz|_\Gamma = 0, \quad
  \Delta\rhoz \in \Linfty , \quad
  \inf\rhoz > 0 , \quad
  \sup\rhoz < 1 ;
  \label{hprhoz}
  \\
  && \xiz \geq 0, \quad
  \sqrt\xiz \in \Huno, \quad
  \lambda(\thetaz,\rhoz) = - \sqrt{\rhoz\xiz}.
  \label{hpxizthetaz}
\Esist
Then, there exist $T>0$ and a triplet $(\rho,\xi,\theta)$
satisfying:
\Bsist
  && \rho \in \H1\Ldue \cap \C0\Huno;
  \label{regrho}
  \\
  && \rho \in \L p{\Wx{2,p}} \;\, 
  \hbox{for each~$p<+\infty$}, \;\,
  \dt\rho \in \LinftyQ,\;Q:=\Omega\!\times\!(0,T); \qquad
  \label{piuregrho}
  \\
  && \xi \in \LinftyQ \cap \W{1,1}\Luno , \quad
  \theta \in \LinftyQ , \quad
  \dt\theta \in \LinftyQ;
  \label{regxitheta}
  \\
  && \inf \rho > 0 , \quad
  \sup\rho < 1 , \quad
  \inf\theta > 0,
  \label{infsup}
\Esist
and solving problem~\pbl.
\Ethm
In the next section, this existence result is proved by a fixed-point argument. 
The method we use seems to be of some interest, because it relies on the 
application of the Tychonoff-Schauder theorem in a weak topology.

\section{Proof of Theorem \ref{Existence}}
\label{Existenceproof}
\setcounter{equation}{0}

By~\eqref{hpcauchy}, we can start by choosing $\epsz>0$ such that
\Beq
  \sqrt{\rhoz\xiz} \leq \cv e^{-1-\cs\rhoz} - 2 \epsz
  \quad \aeO.
  \label{partenza}
\Eeq
Moreover, as the result is local,
we fix a reference final time~$\Tstar>0$ (e.g., $\Tstar=1$)
and assume $T\leq \Tstar$ in the sequel.
Our method is the following. 
By looking at \accorpa{prima}{seconda} and to~\eqref{terza}, separately,
we construct two maps
\Beq
  \FF_1 : \theta \mapsto (\rho,\xi)
  \aand
  \FF_2 : (\rho,\xi) \mapsto \theta
  \non
\Eeq
with proper domains.
Namely, the domain of $\FF_1$
is a convex set $\KK$ depending on~$T$ and on a further parameter~$M$,
and the domain of $\FF_2$ is the range $\RR$ of~$\FF_1$. 
Then, we prove that a suitable choice of $T$ and $M$
ensures that the range of $\FF_2$ is contained in~$\KK$.
This allows us to look for a fixed point of~$\FF_2\circ\FF_1$.
To this aim, we want to use the Tychonoff-Schauder theorem.
For that reason, $\KK$~will be endowed with some weak topology.
We start to construct~$\FF_1$.
The whole argument relies on the technique of~\lavoro\
and the whole paper has to be revisited.
This is done in the next steps.
In particular, as in~\lavoro, we have to consider both the Cauchy problem 
obtained by coupling equation \eqref{seconda} to the second \eqref{cauchy}
and the Cauchy problem
\Beq  
  \dt\xi + \frac{|\dt v|^2 + \source} {\sqrt v} \, \sqrt\xi = 0
  \aand
  \xi(0) = \xiz
  \label{odev}
\Eeq 
(i.e., $\rho$ is replaced by $v$ in \eqref{seconda}),
where $v$ satisfies
\Beq
  v \in D(\Phi) := \graffe{
    v \in \H1\Ldue: \
    v >0 , \
    1/v \in \LinftyQ
  }
  \label{defPhi}
\Eeq
and possibly further conditions later on.
Moreover, the map $\Phi:D(\Phi)\to\L\infty\Luno$
is defined by $\Phi(v)=\sqrt\xi$,
where $\xi$ is the maximal solution to the Cauchy problem~\eqref{odev}.

\Brem
\label{Notation}
In the sequel, our notation is going to reflect dependences, if any, on such parameters as, say, $T$ and $M$; however, possible dependences on problem data such as $\Omega$, $f$, $\source$, and  the initial data, will not be displayed.
\Erem

\step The crucial constants and the maximum principle

We set:
\Beq
  \thetamin := \inf_{0<r<1} \thetam(r)
  \aand
  \thetamax := \sup_{0<r<1} \thetaM(r).
  \label{rangetheta}
\Eeq
A simple calculation yields:
\Beq
  0<\thetamin = e^{-(1+\cs)}
  \aand
  \thetamax = 1;
  \non
\Eeq
thus, we require that $\theta$ obeys the following double limitation: 
\Beq
  \thetamin \leq \theta \leq \thetamax;
  \label{boundstheta}
\Eeq 
in particular, the condition $\inf\theta>0$ (see~\eqref{infsup}) will~automatically hold true.
Now, we notice that \eqref{prima} differs from the analogous one of~\lavoro\
just for the presence of $-\cz\theta$ on the \lhs.
Therefore, even though such a term is space and time dependent,
for a fixed~$\theta$, it can be seen as a part of the smooth perturbation~$\fd$ of the nonlinear term.
Just by thinking of that, the construction of the crucial constants
$\rhomin$, $\rhomax$, and $\ximax$ can be done exactly as in the quoted paper,
provided that the definition of $M_2$ is modified as follows
\Beq
  M_2 := \sup \graffe{ |\fd'(r) - \cz s| : \ \rho\in(0,1), \ s \in(\thetamin,\thetamax) }.
  \label{newM}
\Eeq
Then, the analogues of Lemmas~4.1-4.3 hold in the present case, 
provided that we assume~\eqref{boundstheta}.
Indeed, it suffices to read $\fd'(\rho)-\cz\theta$ in place of $\fd'(\rho)$ in the proofs.

\step The convex set

We define the set $\KK$ as follows:
\Beq
  \KK = \graffe{
    \theta\in\LdueQ: \ \theta,\dt\theta\in\LinftyQ, \ \thetamin\leq\theta\leq\thetamax, \ |\dt\theta|\leq M
  }.
  \label{defK}
\Eeq
Clearly, $\KK$ depends on $T$ and on the a real parameter~$M$,
even though such a dependence is not stressed in the notation.
Moreover, it is non-empty, convex, bounded, and closed.

In the next steps, $\theta\in\KK$ is given 
and we want to solve \accorpa{prima}{seconda} and the first two equations~\eqref{cauchy} for~$(\rho,\xi)$
by applying the procedures of \lavoro.

\step $L^p$ estimates

We make a general observation.
If $p\in(1,+\infty)$ and some function $z$ solves
\Bsist
  && \dt z - \Delta z = g \in L^p(\Qstar)
  \quad \hbox{in $\Qstar:=\Omega\times(0,\Tstar)$}
  \aand \hbox{$\dn z=0$ on the boundary,}
  \non
  \\
  && z(0) = z_0
  \quad \hbox{with} \quad
  z_0 \in \Linfty, \quad
  \Delta z_0 \in \Linfty,
  \aand
  \hbox{$\dn z_0=0$ on the boundary,}
  \non
\Esist
then, the following estimate holds
\Beq
  \norma{\dt z}_{L^p(\Qstar)} + \norma z_{L^p(0,\Tstar;\Wx{2,p})}
  \leq C_p \Bigl( \norma g_{L^p(\Qstar)} + \norma{z_0}_{\Linfty} + \norma{\Delta z_0}_{\Linfty} \Bigr),
  \non
\Eeq
where $C_p$ depends on $\Omega$, $\Tstar$, and~$p$, only.
Therefore, the same constant $C_p$ yields
\Beq
  \norma{\dt z}_{L^p(Q)} + \norma z_{\L p{\Wx{2,p}}}
  \leq C_p \Bigl( \norma g_{L^p(Q)} + \norma{z_0}_{\Linfty} + \norma{\Delta z_0}_{\Linfty} \Bigr)
  \label{genstimap}
\Eeq
for the solution $z$ to the problem
\Beq
  \dt z - \Delta z = g \in L^p(Q)
  \quad \hbox{in $Q$}, \quad
  \hbox{$\dn z=0$ on the boundary}
  \aand
  z(0) = z_0,
  \label{genNeumann}
\Eeq
provided that $T\leq\Tstar$.
Now, note that $z:=\rho$ solves \eqref{genNeumann}
with $g:=-f'(\rho)+\cz\theta+\sqrt{\xi/\rho}$
and that we are assuming $\theta\in\KK$; moreover, \accorpa{hpdati}{hprhoz}~hold.
Thus, by applying~\eqref{genstimap}, we~have:
\Beq
  \norma{\dt\rho}_{L^p(Q)} + \norma\rho_{\L p{\Wx{2,p}}} \leq R_p,
  \label{stimap}
\Eeq
where $R_p$ depends on the same parameters as in \lavoro
and on~$\thetamax$, but not on $M$ or~$T$ (following the rule laid down in Remark~\ref{Notation}, our notation stresses just the dependence on~$p$).
Furthermore, all this is true for the solution $\rho$ coming from
$\sqrt\xi=\Phi(v)$, where $v\in D(\Phi)$ satisfies $\rhomin\leq v\leq\rhomax$.

\step $L^\infty$ estimate

As in~\lavoro, we differentiate with respect to time
and see that $u:=\dt\rho$ satisfies
\Bsist
  && \dt u - \Delta u  {\pier {}+ u} = F + \cz \dt\theta
  \aand \hbox{$\dn u=0$ on the boundary,}
  \label{dtprima}
  \\
  && u(0) = \Delta\rhoz - f'(\rhoz) + \sqrt{\xiz/\rhoz} + \cz\thetaz,
  \label{dtcauchy}
\Esist
where $F$ is as in~\lavoro,
namely,
\Beq
  F := {\pier \dt\rho }- f''(\rho) \dt\rho
  - \frac 12 \, \phi \, \rho^{-3/2} \dt\rho
  - \frac 12 \, \chi \, \bigl( |\dt v|^2 + \source \bigr) (v\rho)^{-1/2}, 
  \label{dtrhs}
\Eeq
where $\phi=\Phi(v)$, with $v\in D(\Phi)$ such that $\rhomin\leq v\leq\rhomax$,
and $\chi$ is some characteristic function.
Then, we fix the right value of $q$ as in~\lavoro,
in order to get the desired estimate.
We have
\Beq
  \norma u_{\LinftyQ}
  \leq C \Bigl( \norma{F+\cz\dt\theta}_{L^q(Q)} + \norma{u(0)}_{\Linfty} \Bigr)
  \non
\Eeq
where $C$ does not depend on~$T$
(by~the above general observation).
As $\theta\in\KK$, we have:
\Beq
  \norma{F+\cz\dt\theta}_{L^q(Q)}
  \leq \norma F_{L^q(Q)} + \cz(|\Omega|T)^{1/q}M;
  \non
\Eeq
by the use of the $L^p$ estimates, we obtain:
\Beq
  \norma{\dt\rho}_{\LinftyQ}
  \leq R_\infty + C_\infty M T^{1/q},
  \label{stimainfty}
\Eeq
where $R_\infty$, $C_\infty$, and $q$ are independent of $M$ and~$T$.

\step The first map

At this point, for every $\theta\in\KK$,
we consider the following problem: in the set
\Bsist
  && \tilde\RR = \Bigl\{ (\rho,\xi) :\
  \rho \in \H1\Ldue \cap \C0\Huno , \
  \xi \in \W{1,1}\Luno
  \hskip3em
  \non
  \\
  && \hphantom{\KK' = \Bigl\{}
  \rhomin \leq \rho \leq \rhomax , \
  0 \leq \xi \leq \ximax
  \non
  \\
  && \hphantom{\KK' = \Bigl\{}
  \norma{\dt\rho}_{L^p(Q)} + \norma\rho_{\L p{\Wx{2,p}}} \leq R_p \
  \hbox{for every $p\in(1,+\infty)$}
  \non
  \\
  && \hphantom{\KK' = \Bigl\{}
  \norma{\dt\rho}_{\LinftyQ}
  \leq R_\infty + C_\infty M T^{1/q}
  \Bigr\},
  \label{defRtilde}
\Esist
find a pair 
\Beq
  \hbox{$(\rho,\xi)$ 
    satisfying \eqref{prima}, \eqref{seconda},
    the first two of \eqref{cauchy}, and \eqref{maximality}}
  \label{predefFuno}
\Eeq
Indeed, the analogue of the map $\Psi$ considered in \lavoro\
can be defined in the same way and actually maps its domain into itself
also in the present case, because of our choice of the constants.
Moreover, in performing the contraction estimate,
just differences have to be considered and $\theta$ disappears.
Therefore, by setting:
\Beq
  \hbox{for $\theta\in\KK$}, \quad
  \hbox{$\FF_1(\theta)$ is the solution $(\rho,\xi)\in\tilde\RR$ \
    to problem \eqref{predefFuno},}
  \label{defFuno}
\Eeq
we obtain a well-defined map $\FF_1:\KK\to\tilde\RR$.
We~set:
\Beq
  \hbox{$\RR$ is the range of $\FF_1$}.
  \label{defR}
\Eeq
Such a subset of $\tilde\RR$ depends on $T$ and $M$ and
our next step is the definition of a map $\FF_2:\RR\to\KK$
that works as follows:
given $(\rho,\xi)\in\RR$, the value $\FF_2(\rho,\xi)$
is a function $\theta$ satisfying~\eqref{terza} and the third of~\eqref{cauchy}.
However, before doing that, some more work is needed,
since we are not sure that such a map is well-defined.
It is not, indeed, unless $T$ and $M$ satisfy suitable constraints,
and the next steps are devoted to find~them.

\step Higher regularity estimates

Let $(\rho,\xi)\in\RR$.
Then, $(\rho,\xi)$ solves~\eqref{predefFuno}
for some $\theta\in\KK$, whence $u:=\dt\rho$
solves \accorpa{dtprima}{dtcauchy}
where $F$ is the same as~\eqref{dtrhs}.
On noting that \accorpa{hprhoz}{hpxizthetaz} imply
$u(0)\in\Huno$, we have by the above $L^p$ estimates
and the general theory that
\Bsist
  && \norma{\dt^2\rho}_{\LdueQ}
  + \norma{\dt\rho}_{\L2\Huno}
  = \norma{\dt u}_{\LdueQ}
  + \norma u_{\L2\Huno}
  \non
  \\
  && \leq c \bigl(
    \norma{F+\cz\dt\theta}_{\LdueQ}
    + \norma{u(0)}_{\Huno}
  \bigr)
  \leq C'(T,M),
  \label{higherreg}
\Esist
where the form of the dependence of $C'(T,M)$ on $T$ and~$M$ 
is not important in the sequel.

\step New a priori estimate

Let $(\rho,\xi)\in\RR$.
Then, $(\rho,\xi)$ solves~\eqref{predefFuno}
for some $\theta\in\KK$.
Therefore, as in \lavoro, we~have
\Beq
  \dt\sqrt\xi = - \chi \, \frac{|\dt\rho|^2 + \source} {2\sqrt\rho},
  \non
\Eeq
where $\chi$ is some characteristic function, whence immediately
\Beq
  \norma{\dt\sqrt\xi}_{\LinftyQ}
  \leq \frac 1 {2\sqrt\rhomin} \,
    \bigl( \norma{\dt\rho}_{\LinftyQ}^2 + \norma\source_{L^\infty(\Omega\times(0,+\infty))}
    \bigr).
  \non
\Eeq
Hence, if $c$ stands for different constants independent of $M$ and~$T$,
we deduce:
\Bsist
  && \modulo{\dt\sqrhoxi}
  = \modulo{ \frac {\dt\rho} {2\sqrt\rho} \, \sqrt\xi
    + \sqrt\rho \, \dt\sqrt\xi
  }
  \leq \frac {\sqrt\ximax} {2\sqrt\rhomin} \, \norma{\dt\rho}_{\LinftyQ}
  + \sqrt\rhomax \norma{\dt\sqrt\xi}_{\LinftyQ}
  \qquad
  \non
  \\
  && \leq c \Bigl( 1 + \norma{\dt\rho}_{\LinftyQ}^2 \Bigr)
  \leq c \Bigl( 1 + \bigl( R_\infty + C_\infty M T^{1/q} \bigr)^2 \Bigr),
  \non
\Esist
and we conclude that
\Beq
  \modulo{\dt\sqrhoxi} \leq C_1 \bigl( 1 + M^2 T^{2/q} \bigr),  \label{stimadtsqrt}
\Eeq
where $C_1$ is independent of $T$ and~$M$.

\step First restriction on parameters

Observe that
\Beq
  \modulo{ \dt \bigl( \sqrhoxi - \cv e^{-1-\cs\rho} \bigr) }
  \leq |\dt\sqrhoxi| + \cv \cs |\dt\rho|.
  \non
\Eeq
By accounting for \eqref{stimadtsqrt} and~\eqref{stimainfty}, we deduce~that
\Beq
  \norma{ \dt \bigl( \sqrhoxi - \cv e^{-1-\cs\rho} \bigr)}_{\LinftyQ}
  \leq \frac {C_2} 2 \, \bigl( 1 + M^2 T^{2/q} + M T^{1/q} \bigr)
  \leq C_2 \bigl( 1 + M^2 T^{2/q} \bigr),
  \non
\Eeq
where $C_2$ is independent of $T$ and~$M$;
this implies that
\Beq
  \norma{
    \bigl( \sqrhoxi - \cv e^{-1-\cs\rho} \bigr)
    - \bigl( \sqrt{\rhoz\xiz} - \cv e^{-1-\cs\rhoz} \bigr)
  }_{\LinftyQ}
  \leq C_2 T \bigl( 1 + M^2 T^{2/q} \bigr).
  \non
\Eeq
By~\eqref{partenza},
we conclude that
\Beq
  \sqrhoxi \leq \cv e^{-1-\cs\rho} - \epsz
  \quad \aeQ,
  \label{stimaepsz}
\Eeq
whenever $M$ and $T$ satisfy
\Beq
  C_2 T \bigl( 1 + M^2 T^{2/q} \bigr) \leq \epsz \,.
  \label{primarestriz}
\Eeq
Therefore, if $T$ and $M$ satisfy~\eqref{primarestriz},
every $(\rho,\xi)\in\RR$ fulfils~\eqref{stimaepsz}.

\step The second map

Assume $(\rho,\theta)\in\RR$ and \eqref{primarestriz}.
By~\eqref{stimaepsz}, we have, in particular, that
\Beq
  \sqrhoxi < \cv e^{-1-\cs\rho}
  \quad \aeQ.
  \non
\Eeq
Therefore, \aaQ, the equation
$\lambda(\rho(x,t),s)=-\sqrt{\rho(x,t)\xi(x,t)}$
has two solution $s_1$  and $s_2$
which belong to $(0,\thetam(\rho(x,t)))$
and $(\thetam(\rho(x,t)),\thetaM(\rho(x,t)))$, respectively.
We term the latter $\theta_2(x,t)$ for a while
and obtain a bounded function,
thus a function $\theta_2\in\LdueQ$.
Then, we can define $\FF_2:\RR\to\LdueQ$ by setting:
$\FF_2(\rho,\xi)$ is such a function~$\theta_2$.
Therefore, for $(\rho,\xi)\in\RR$,
\Beq
  \theta = \FF_2(\rho,\xi)
  \quad \hbox{means} \quad
  \thetam(\rho) < \theta < \thetaM(\rho)
  \aand
  \lambda(\rho,\theta) = - \sqrhoxi
  \quad \aeQ .
  \label{defFdue}
\Eeq
As both $\rho$ and $\xi$ are continuous with respect to time
(\aaO) and the function $\theta:=\FF_2(\rho,\xi)$ 
is always different from~$\thetam(\rho)$, 
time continuity holds for $\theta$ as well,
and it is clear that $\theta(0)=\thetaz$.

Although $\FF_2$ is well-defined,
we still have to find the restriction on $T$ and $M$
that ensures that the range of $\FF_2$ is contained in~$\KK$.

\step Estimate from below

Let $(\rho,\xi)\in\RR$, set $\theta:=\FF_2(\rho,\xi)$,
and assume~\eqref{primarestriz}.
Then \eqref{stimaepsz} holds.
For almost every $(x,t)\in Q$,
we~write the second--order Taylor expansion 
of the function $s\mapsto\lambda(\rho(x,t),s)$
with center at $s=\thetam(\rho(x,t))$ (hereafter, to lighten our notation, we refrain from writing $(x,t)$).
We find $s'\in(\thetam(\rho),\theta)$ such that the following holds:
\Bsist
  && - \sqrhoxi
  = \lambda(\rho,\theta)
  = \lambda(\rho,\thetam(\rho))
  + \frac{\partial\lambda}{\partial s} (\rho,\thetam(\rho)) \, (\theta - \thetam(\rho))
  + \frac 12 \, \frac{\partial^2\lambda}{\partial s^2} (\rho,s') \, (\theta - \thetam(\rho))^2
  \non
  \\
  && = - \cv e^{-1-\cs\rho} + \frac\cv {2s'} \, (\theta - \thetam(\rho))^2,
  \non
\Esist
in view of \eqref{defthetamM} and~\eqref{thetamM}.
By~\eqref{primarestriz}, we~deduce that
\[
  \frac\cv {2s'} \, (\theta - \thetam(\rho))^2
  = \cv e^{-1-\cs\rho} - \sqrhoxi
  \geq \epsz,
\]
whence
\[
(\theta - \thetam(\rho))^2
  \geq \frac {2 \epsz} \cv \, s'
  \geq \frac {2 \epsz} \cv \, \thetam(\rho)
  = \frac {2\epsz} \cv \, e^{-1-\cs\rho}
  \geq \frac {2\epsz} \cv \, e^{-1-\cs\rhomax} .
\]
As $\theta>\thetam(\rho)$, we conclude that
\Beq
  \theta - \thetam(\rho)
  \geq 2\deltaz
  \quad \aeQ
  \label{dalbasso}
\Eeq
(with an obvious definition of~$\deltaz$; note that $\deltaz$ is independent of $T$ and~$M$).

\Blem
\label{Stimadsl}
Assume that $r\in(\rhomin,\rhomax)$, $s\in(\thetamin,\thetamax)$, and $\delta>0$.
Then,
\Beq
  s \geq \thetam(r) + \delta
  \quad \hbox{implies} \quad
  \frac {\partial\lambda} {\partial s} (r,s) 
  \geq \cv \ln(1+\delta/\thetamax) \,.
  \label{stimadsl}
\Eeq
\Elem

\Bdim
Indeed, if $s\geq\thetam(r)+\delta$, we have that
\Bsist
  && \frac {\partial\lambda} {\partial s} (r,s) 
  = \cz r + \cv + \cv \ln s
  = \cz r + \cv + \cv \ln\thetam(r) + \cv \ln \frac s {\thetam(r)}
  \non
  \\
  && =  \frac {\partial\lambda} {\partial s} (r,\thetam(r))
  + \cv \ln \frac s {\thetam(r)}
  = \cv \ln \frac s {\thetam(r)}
  \geq \cv \ln \frac {\thetam(r) + \delta} {\thetam(r)}
  = \cv \ln \bigl( 1 + \delta/\thetam(r) \bigr)
  \non
\Esist
and the desired conclusion follows, because $\thetam(r)\leq\thetamax$.
\Edim

\step Estimate of the time derivative

Assume that $(\rho,\xi)\in\RR$, $\theta:=\FF_2(\rho,\xi)$,
and that~\eqref{primarestriz} holds.
On the one hand, we~have that
\Beq
  \frac {\partial\lambda} {\partial s} (\rho,\theta) \, \dt\theta
  = - \dt \sqrhoxi
  - \frac{\partial\lambda}{\partial r} (\rho,\theta) \, \dt\rho
  = - \dt \sqrhoxi
  - \cz \theta \dt\rho
  \quad \aeQ .
  \non
\Eeq
On the other hand, \eqref{dalbasso}~holds.
Therefore, if we apply the previous lemma for $\delta:=2\deltaz$ 
and set $\Lz:=\bigl(\cv\ln(1+2\deltaz/\thetamax)\bigr)^{-1}$, we~obtain:
\Beq
  |\dt\theta|
  \leq \Lz \bigl( |\dt\sqrhoxi| + \cz \theta |\dt\rho| \bigr)
  \leq \Lz \bigl( |\dt\sqrhoxi| + \cz \thetamax |\dt\rho| \bigr) .
  \label{perstimadttheta}
\Eeq
Now, we account for~\eqref{stimadtsqrt} and conclude that
\Beq
  |\dt\theta| 
  \leq \Lz C_1 (1 + M^2 T^{2/q})
  + \Lz \cz \thetamax (R_\infty + C_\infty M T^{1/q})
  \quad \aeQ.
  \label{stimadttheta}
\Eeq
At this point, we are ready to choose the constants $T$ and~$M$.

\step Choice of the constants

Clearly, \eqref{stimadttheta} implies $|\dt\theta|\leq M$, whenever
\Beq
  \Lz C_1 (1 + M^2 T^{2/q}) 
  + \Lz \cz \thetamax (R_\infty + C_\infty M T^{1/q})
  \leq M .
  \label{secondarestriz}
\Eeq
Therefore, we choose, e.g., $M=\Lz C_1+\Lz\cz\thetamax R_\infty+1$ and $T$
small enough for both 
\eqref{primarestriz} and \eqref{secondarestriz} to hold.
With such a choice, 
the range of $\FF_1$ is contained in~$\RR$ and 
$\FF_2$ is a well-defined map whose range is contained in~$\KK$.
Hence, $\FF:=\FF_2\circ\FF_1$ is a well-defined map  from $\KK$ into $\LdueQ$
that satisfies $\FF(\KK)\subseteq\KK$.

\step Conclusion of the proof
 
As anticipated, we aim to an application of the Tychonoff-Schauder fixed point theorem.
As far as the topology of $\KK$ is concerned, we see $\KK$ as a 
subset of the topological vector space obtained by endowing $\LdueQ$ with the weak topology.
Therefore, the convex set $\KK$ is compact.
So, the last point of our proof is the continuity of $\FF$
with respect to the topology of~$\KK$.
To this aim, we observe that
$\LdueQ$ with its strong topology is reflexive and separable and that $\KK$ is bounded.
Thus, the topology of $\KK$ comes from a metric.
In particular, $\FF$ is continuous
if and only if it is sequentially continuous.
So, we pick $\hthetan,\htheta$ such that
\Beq
  \hthetan \in \KK
  \quad \hbox{for every $n$}
  \aand
  \hthetan \to \htheta
  \quad \hbox{weakly in $\LdueQ$}
  \label{hpconv}
\Eeq
and pass to prove that the corresponding functions
$\thetan:=\FF(\hthetan)$ and $\theta:=\FF(\htheta)$ satisfy
\Beq
  \thetan \to \theta
  \quad \hbox{weakly in $\LdueQ$}.
  \label{tesi}
\Eeq
It suffices for us to show that every subsequence of $\graffe{\thetan}$
has a subsequence that converges to $\theta$ weakly in~$\LdueQ$.
In the sequel, to lighten our notation,
we denote by the same symbol ($\graffe{y_n}$, say) both a sequence and all its subsequences.
So, we start by any subsequence $\graffe{\thetan}$ of the given sequence
and look for a subsequence of it satisfying~\eqref{tesi}.
We set for convenience $(\rhon,\xin):=\FF_1(\hthetan)$.
Therefore, we have $(\rhon,\xin)\in\RR\subset\tilde\RR$
and all the estimates given in definition~\eqref{defRtilde} hold true,
as well as~\eqref{higherreg}. 
In particular, we have
\Bsist
  && \rhon \to \rho
  \quad \hbox{weakly in $\W{1,p}{\Lx p}\cap\L p{\Wx{2,p}}$ for every $p<+\infty$}
  \qquad\quad
  \label{convrhon}
  \\
  && \dt\rhon \to \dt\rho
  \quad \hbox{weakly in $\H1\Ldue\cap\L2\Huno$}
  \label{convdtrhon}
\Esist
for some $\rho$, at least for a subsequence.
From \accorpa{convrhon}{convdtrhon} we also~have that
\Beq
  \rhon \to \rho
  \quad \hbox{uniformly in $Q$}
  \aand
  \dt\rhon \to \dt\rho
  \quad \hbox{strongly in $\LunoQ$}.
  \label{rhonok}
\Eeq
As $\rhomin\leq\rhon\leq\rhomax$ for every~$n$,
we infer the uniform convergence of $f'(\rhon)$ to~$f'(\rho)$.
Now, we consider the maximal solution $\xi$ to the Cauchy problem
obtained by coupling \eqref{seconda} (with such a~$\rho$)
and the second~\eqref{cauchy}.
We prove that $(\rho,\xi)=\FF_1(\htheta)$.
By using \cite[Lemma~4.8]{CGPS},
we easily deduce that, \aeO\ and for every $n$ and $t\in[0,T]$,
\Bsist
  && |\sqrt{\xin(t)} - \sqrt{\xi(t)}|
  \leq c \ioT \bigl| \, |\dt\rhon(s)|^2 - |\dt\rho(s)|^2 \bigr| \, ds
  + c \ioT |\rhon(s) - \rho(s)| \, ds
  \non
  \\
  && \leq c \ioT |\dt\rhon(s) - \dt\rho(s)| \, ds
  + c \ioT |\rhon(s) - \rho(s)| \, ds,
  \non
\Esist
where $c$ stands for different constants independent of~$n$.
We infer~that
\Beq
  \norma{\sqrt{\xin}-\sqrt{\xi}}_{\L\infty\Luno}
  \leq c \bigl( 
    \norma{\dt\rhon-\dt\rho}_{\LunoQ}
    + \norma{\rhon-\rho}_{\LunoQ},
  \bigr)
  \non
\Eeq
and, owing to~\eqref{rhonok}, we deduce that 
$\sqrt\xin$ converges to $\sqrt\xi$ in~$\L\infty\Luno$,
thus \aeQ\ at least for a subsequence.
At this point, it is clear that
$(\rho,\xi)$ satisfies problem \pbl\ where we read $\htheta$ in place of~$\theta$.
Indeed, \eqref{seconda}, the second~\eqref{cauchy}, and \eqref{maximality}
hold by definition of~$\xi$.
On the other hand, it is \sfw\ to let $n$ tend to infinity in the variational equation,
or in an integrated version of it, and in the Cauchy conditions.
We conclude that $(\rho,\xi)=\FF_1(\htheta)$,
and the next step is the convergence of~$\graffe{\thetan}$ to $\theta$ weakly in~$\LdueQ$,
at least for a subsequence. As a matter of fact, we can prove a strong--convergence result.
Recall that the estimate from below~\eqref{dalbasso}
holds for~$\thetan$.
It follows~that
\Beq
  \thetan \geq \thetam(\rhon) + 2\deltaz
  = e^{-1-\cs\rhon} + 2\deltaz
  = e^{-1-\cs\rhon} - e^{-1-\cs\rho} + \thetam(\rho) + 2\deltaz,
  \non
\Eeq
whence
\Beq
  \thetan \geq \thetam(\rho) + \deltaz
  \quad \hbox{provided that} \quad
  |\rhon - \rho| \leq \etaz,
  \non
\Eeq
where $\etaz>0$ is given by the uniform continuity 
of the exponential on every bounded interval
(namely, $r,r'\in[\rhomin,\rhomax]$ and $|r-r'|\leq\etaz$
imply that $|e^{-1-\cs r}-e^{-1-\cs r'}|\leq\deltaz$).
Moreover, the similar inequality 
\Beq
  s \geq \thetam(r) + \deltaz
  \non
\Eeq
holds true whenever $|\rhon-\rho|\leq\etaz$ and
\Beq
  \min\graffe{\rhon,\rho} \leq r \leq \max\graffe{\rhon,\rho}
  \aand
  \min\graffe{\thetan,\theta} \leq s \leq \max\graffe{\thetan,\theta}.
  \label{framinmax}
\Eeq
On the other hand, the uniform convergence~\eqref{rhonok}
implies that $|\rhon-\rho|\leq\etaz$ in~$Q$ for $n$ large enough.
Therefore, for such values of~$n$, the following holds.
We apply the Lagrange mean value theorem to~$\lambda$.
For a.a.\ $(x,t)\in Q$ (once again, we omit writing $(x,t)$ in the sequel)
and~suitable (space and time dependent) $r_n$~and $s_n$ satisfying~\eqref{framinmax},
we~have that
\Beq
  \sqrhonxin - \sqrhoxi
  = \lambda(\rho,\theta) - \lambda(\rhon,\thetan)
  = \frac {\partial\lambda} {\partial r} (r_n,s_n) \, (\rho - \rhon)
  + \frac {\partial\lambda} {\partial s} (r_n,s_n) \, (\theta - \thetan).
  \non
\Eeq
Now, we observe that $r_n\in[\rhomin,\rhomax]$ and $s_n\geq\thetam(r_n)+\deltaz$
since $|\rho-r_n|\leq\etaz$.
Hence, we can apply Lemma~\ref{Stimadsl} and arrive at \eqref{stimadsl}
with $r=r_n$, $s=s_n$, and $\delta=\deltaz$.
This implies that
\Beq
  \ln(1 + \deltaz/\thetamax) \, |\thetan - \theta|
  \leq |\sqrhonxin - \sqrhoxi| 
  + \cz |\rhon - \rho|.
  \non 
\Eeq
As the \rhs\ tends to $0$ \aeQ\ and is bounded in~$\LinftyQ$,
we deduce that $\thetan$ strongly converges  to~$\theta$ in~$\LdueQ$.
This concludes the proof.



\Begin{thebibliography}{10}

\bibitem{CN} B.D. Coleman and W. Noll,  The thermodynamics of elastic 
materials with heat conduction and viscosity. Arch. Rational Mech. Anal. 
\textbf{13} (1963) 167--178. 

\bibitem{CDKMS}
P. Colli, A. Damlamian, N. Kenmochi, M. Mimura, J. Sprekels (eds.), 
``Proceedings of International Conference on: Nonlinear Phenomena with Energy 
Dissipation. Mathematical Analysis, Modeling and Simulation'',  
GAKUTO International Series. Mathematical Sciences and Applications~{\bf 29}, {Gakk\=otosho}, Tokyo 2008. 

\bibitem{CGPS} 
P. Colli, G. Gilardi, P. Podio-Guidugli, J. Sprekels, {\em
Existence and uniqueness of a global-in-time solution
to a phase segregation problem of the Allen-Cahn type}, 
Math. Models Methods Appl. Sci. {\pier {\bf 20} (2010)
519--541.}

\bibitem{CKS}
P. Colli, N. Kenmochi \& J. Sprekels (eds.), ``Dissipative 
Phase Transitions'',  Series on Advances in Mathematics for Applied 
Sciences~{\bf 71}, World Scientific Publishing Co. Pte. Ltd., Hackensack, 
NJ 2006. 

\bibitem{fremond}
M. Fr\'emond,
``Non-smooth Thermomechanics'',
Springer-Verlag, Berlin, 2002.

\bibitem{gurtin}
M.E. Gurtin,
{\em Generalized Ginzburg-Landau and 
Cahn-Hilliard equations based on a microforce balance},
Phys.~D
{\bf 92} (1996)
178--192.

\bibitem{miranv} 
A. Miranville,
{\em Consistent models of {C}ahn-{H}illiard-{G}urtin equations with
{N}eumann boundary conditions},
Phys. D {\bf 158} (2001) 233--257.

\bibitem{Podio}
P. Podio-Guidugli, 
{\em Models of phase segregation and diffusion of atomic species on a lattice},
Ric. Mat. {\bf 55} (2006) 105--118.

\End{thebibliography}

\end{document}